\documentclass[a4paper,10pt]{article}

\usepackage{amsmath}
\usepackage{graphicx}
\usepackage{latexsym}
\usepackage{amssymb}
\usepackage{amsbsy}

\DeclareSymbolFont{msbm}{U}{msb}{m}{n}
\DeclareMathSymbol{\C}{\mathalpha}{msbm}{'103}
\DeclareMathSymbol{\R}{\mathalpha}{msbm}{'122}
\DeclareMathSymbol{\Z}{\mathalpha}{msbm}{'132}
\DeclareMathSymbol{\N}{\mathalpha}{msbm}{'116}

\newtheorem{remark}{Remark}

\newtheorem{theorem}{Theorem}
\newtheorem{definition}{Definition}
\newtheorem{lemma}{Lemma}
\newtheorem{corollary}{Corollary}

\def\RR{\mathbb R}

\def\be{\begin{equation}}
\def\ee{\end{equation}}
\def\bea{\begin{eqnarray}}
\def\ba{\begin{array}{l}\displaystyle}
\def\eea{\end{eqnarray}}
\def\ea{\end{array}}

\newcommand{\tA}{\widetilde{A}}
\newcommand{\ta}{\widetilde{a}}
\newcommand{\tw}{\widetilde{w}}
\newcommand{\tc}{\widetilde{c}}

\begin{document}
\title{High order asymptotic-preserving schemes for the Boltzmann equation}

\author{Giacomo Dimarco\thanks{Universit\'{e} de Toulouse; UPS, INSA, UT1, UTM;
CNRS, UMR 5219; Institut de Math\'{e}matiques de Toulouse; F-31062
Toulouse, France. ({\tt giacomo.dimarco@math.univ-toulouse.fr}).}
\and Lorenzo Pareschi\thanks{Mathematics Department, University of
Ferrara and CMCS, Ferrara, Italy ({\tt lorenzo.pareschi@unife.it}).}
}
\date{}
\maketitle

\begin{abstract}
\noindent In this note we discuss the construction of high order
asymptotic preserving numerical schemes for the Boltzmann
equation. The methods are based on the use of Implicit-Explicit
(IMEX) Runge-Kutta methods combined with a penalization technique
recently introduced in~\cite{Filbet}.
\medskip


\begin{center}
{\large {\bf Sch\'emas d'ordre \'el\'ev\'e et pr\'eservant
l'asymptotique pour l'\'equation de Boltzmann}}
\end{center}
\begin{center}
{\bf R\'esum\'e}
\end{center}

\noindent Dans cette note nous discutons la construction de
sch\'{e}mas d'ordre \'{e}lev\'{e} pour l'\'{e}quation de Boltzmann
qui pr\'{e}servent la limite asymptotique. Les m\'{e}thodes sont
bas\'{e}es sur l'utilisation de sch\'{e}mas de Runge-Kutta
explicites-implicites combin\'{e}es avec une technique de
p\'{e}nalisation introduit r\'{e}cemment par~\cite{Filbet}.
\medskip

\noindent {\bf Mots-cl\'{e}s}~: M\'{e}thodes Runge-Kutta
Implicites-Explicites, \'{e}quations raides, \'equation de
Boltzmann, limite fluide, sch\'{e}mas pr\'{e}servant l'asymptotique.
\end{abstract}

%\tableofcontents
\bigskip

{\Large {\bf Version fran\c{c}aise abr\'eg\'ee}}
\bigskip

\noindent Les \'{e}quations cinétiques, comme l'\'{e}quation de
Boltzmann sont utilis\'{e} avec succ\`{e}s dans de nombreuses
applications r\'{e}elles. L'\'{e}quation de Boltzmann d\'{e}crit
l'\'{e}volution temporelle de la fonction de distribution d'un gaz
avec des interactions binaires \'{e}lastiques. Il est importante de
mentionner que la solution num\'{e}rique de l'op\'{e}rateur de
collision repr\'{e}sente un d\'{e}fi majeur pour les m\'{e}thodes
num\'{e}riques traditionnelles qui n'est pas encore r\'{e}solu. Cela
est particuli\`{e}rement vrai en proximit\'{e} des r\'{e}gimes
fluides. Dans ces r\'{e}gimes le taux des collisions
intermol\'{e}culaires croît de façon exponentielle et donc le temps
entre deux collisions successives devient tr\`{e}s petit. D'autre
part, l'\'{e}chelle de temps r\'{e}el pour l'\'{e}volution du gaz
est l'\'{e}chelle de temps de la dynamique des fluides, qui est
normalement beaucoup plus grande que le temps entre deux collisions.
Une mesure de l'importance des collisions est donn\'{e}e par le
nombre de Knudsen $\varepsilon $, qui est grand dans la limite
rar\'{e}fi\'{e}e et petit dans la limite fluide. Ainsi, les
approches num\'{e}riques traditionnelles perdent leur efficacit\'{e}
en raison de la n\'{e}cessit\'{e} d'utiliser de temps tr\`{e}s
petits pour la discretization temporelle. Nous rappelons que la
discr\'{e}tisation directe en temps de l'\'{e}quation de Boltzmann
est un gros probl\`{e}me dans les r\'{e}gimes raides en raison de la
haute dimensionnalit\'{e} et de la non-lin\'{e}arit\'{e} de
l'op\'{e}rateur de collision qui rend peu pratique l'utilisation de
solveurs implicites.

Plusieurs auteurs ont abord\'{e} le probl\`{e}me dans le r\'{e}cent
pass\'{e}( voir~\cite{dimarco6, DiPa12, Filbet, Jin2, Lem, PRimex})
et les r\'{e}f\'{e}rences \'{a} l'int\'{e}rieur). Une strat\'{e}gie,
parmi les plus puissantes, consiste en la construction des
sch\'{e}mas dits pr\'{e}servant l'asymptotique. Ces techniques
permettent de r\'{e}soudre le probl\`{e}me dans tout le domaine pour
tous les choix de pas de temps et de nombre de Knudsen. Dans cette
note, nous proposons une nouvelle classe de schemas Runge-Kutta
Implicites-Explicites pour l'\'{e}quation de Boltzmann. Pour
construire nos schemas, nous consid\'{e}rons une d\'{e}composition
du terme de gain de l'op\'{e}rateur de collision en un partie en
\'{e}quilibre et en une partie en non \'{e}quilibre. Cette
d\'{e}composition de l'int\'{e}grale de Boltzmann a \'{e}t\'{e}
\'{e}galement introduite par Jin et Filbet dans~\cite{Filbet}. Les
principaux avantages de l'approche propos\'{e}e ici est que cela
fonctionne de mani\`{e}re uniforme pour une large gamme de nombres
de Knudsen et \'{e}vite la solution d'un syst\`{e}me d'\'{e}quations
non lin\'{e}aires, m\^{e}me dans les r\'{e}gimes raides. De m\^{e}me
que pour~\cite{dimarco6}, nous obtenons des conditions suffisantes
pour la stabilit\'{e} asymptotique et la pr\'{e}servation
asymptotique de l'ordre temporel des schemas. En plus, nous
construirons les sch\'{e}mas tels qu'ils pr\'{e}servent la
positivit\'{e} des solutions et les quantit\'{e}s physiques
conserv\'{e}es. Pour plus de d\'etails nous renvoyons
\`a~\cite{DiPa12}.

%%%%%%%%%%%%%%%%%%%%%%%%%%%%%%%%%%%%%%%%%%%%%%%%%%%%%%%%%%%%%%%%%%%%%%%%%%%%%%%%%%%%%%%%%%%%%%
\section{Introduction}
The computation of fluid-kinetic interfaces and asymptotic behaviors
involves multiple scales where most numerical methods lose their
efficiency because they are forced to operate on a very short time
scale (see~\cite{dimarco6, DiPa12, Filbet, Jin2, Lem, PRimex} and the references therein
for a more complete bibliography). The Boltzmann equation close to
fluid regimes represents the prototype example \be
\partial_t f + v\cdot\nabla_{x}f
=\frac{1}{\varepsilon}Q(f,f).\label{eq:1}\ee Here $f(x,v,t)$ is a
non negative function describing the time evolution of the
distribution of particles with velocity $v \in \R^{3}$ and
position $x \in \Omega \subset \R^{d_x}$ at time $ t
> 0$.

The operator $Q(f,f)$ describes the particles interactions. In the
general case of the Boltzmann binary collision, it has the form
\be Q_{B}(f,f)=\int_{\RR^3\times S^2} B(|v-v_*|,n)
[f(v')f(v'_*)-f(v)f(v_*)]\,dv_*\,dn \ee where \be
v'=v+\frac12(v-v_*)+\frac12|v-v_*|n,\quad
v'_*=v+\frac12(v-v_*)-\frac12|v-v_*|n, \ee and $B(|v-v_*|,n)$ is a
nonnegative collision kernel characterizing the details of the
collision.

The Knudsen number $\varepsilon>0$ is a non dimensional measure of
the importance of collisions and is large in rarefied regions and
small where the system is close to the fluid limit. In the latter
regime, the intermolecular collision rate grows quickly and thus
the collisional time scale becomes very small. On the other hand,
the actual time scale of the  evolution is the fluid dynamic
scale, which can be much larger than the collisional time.

In fact, for small values of $\varepsilon$ the distribution
function is well approximated by a local Maxwellian \be
M[f]=\frac{\rho}{(2\pi
T)^{3/2}}\exp\left(\frac{-|w-v|^{2}}{2T}\right), \label{eq:M}\ee
where $\rho$, $w$, $T$ are the density, mean velocity and
temperature of the gas in the x-position and at time $t$ defined
as \be (\rho,\rho w,E)^T=\int_{\RR^3} f
\left(1,v,\frac{v^2}{2}\right)^T\,dv, \qquad
T=\frac1{3\rho}(E-\rho|w|^2). \ee

Now, passing to the limit for $\varepsilon\rightarrow 0$ and
integrating (\ref{eq:1}) against $1$, $v$ and $v^2$ we recover the
system of compressible Euler equations \be
\partial_t u+\nabla_x\cdot F(u)=0
\label{eq:Euler} \ee with
\[
u=(\rho,w,E)^T,\qquad F(u)=(\rho w, \varrho w \otimes (w+pI),
Ew+pw)^T,\quad p=\rho T,
\]
where $I$ is the identity matrix.

Implicit-Explicit (\emph{IMEX}) Runge-Kutta schemes represent a
powerful tool for the numerical treatment of stiff terms in PDEs~\cite{dimarco6, boscarino, PRimex}. When necessary they can be designed in
order to achieve suitable asymptotic preserving (\emph{AP})
properties. Their direct application to the Boltzmann equation
however is not trivial since the complicated nonlinear structure of
the collisional operator makes prohibitively expensive the use of
implicit solvers for the stiff collision term. Additional
difficulties are given by the need to preserve the most relevant
physical properties of the solution, like conservation of mass,
momentum and energy, nonnegativity and entropy inequality. In this
short note we will illustrate how the introduction of a suitable
penalization technique as in~\cite{Filbet} permits to extend
succesfully the IMEX formalism also to the challenging case of the
Boltzmann equation.
%Let us finally mention that several authors have
%tackled the above problem in the recent past, we refer to
%\cite{dimarco6, Lem} and the references therein.

\section{IMEX schemes for the Boltzmann equation}
In order to apply efficiently the IMEX Runge-Kutta approach to the Boltzmann equation we must avoid the prohibitive cost of the implicit evaluation of the stiff collision term. In order to achieve this we first reformulate the collision part using a suitable penalization term.

\subsection{Decomposition of the collision integral}  First, we observe that we
can rewrite $Q_{B}(f,f)$ as~\cite{toscani} \be
Q_{B}(f,f)=\frac{1}{\varepsilon}(P(f,f)-\mu f),\label{eq:10} \ee
where $P(f,f)=Q_{B}(f,f)+\mu f$ and $\mu>0$ is a constant such that
$P(f,f)\geq 0$.
%Typically $\mu f$ is an upper bound of the loss part of the collision term  \be \mu\geq \int_{\RR^3\times S^2} B(|v-v_*|,n)
%f(v_*)\,dv_*\,dn. \label{eq:P}\ee

Observe that, by construction, the following property is
verified by the operator $P(f,f)$\be \frac1{\mu}\int_{\R^3}
P(f,f)\left(1,v,\frac{v^2}{2}\right)^T\,dv = \int_{\R^3}
f\left(1,v,\frac{v^2}{2}\right)^T\,dv=u. \ee Thus, $P(f,f)/\mu$ is a
density function and we can consider the following decomposition \be
P(f,f)/\mu=M[f]+g,\label{eq:dec}\ee where the function $g$
represents the deviations from equilibrium of $P(f,f)$.
% and
%from the definition above it follows that $g$ is in general non
%positive and with vanishing first moments.

Thus the collision operator can be rewritten in the form
\be
Q_{B}(f,f)=\frac{\mu}{\varepsilon}g+\frac{\mu}{\varepsilon}(M[f]-f)=\frac{\mu}{\varepsilon}\left(\frac{P(f,f)}{\mu}-M[f]\right)+
\frac{\mu}{\varepsilon}(M[f]-f).\label{eq:11}\ee The above
reformulation is equivalent to the penalization method for the
collision operator recently introduced in~\cite{Filbet}. Clearly, since the
problem is stiff as a whole a fully implicit method should be
used in the numerical integration to avoid stability constraints of
the type $\Delta t = O(\varepsilon)$. On the other hand, the linear
part itself $(M[f]-f)$ suffices to characterize the correct large time behavior
of $f$. Therefore, instead of fully implicit methods, one may use
methods which are implicit in the linear part and explicit in the
non-linear part. This however, as we will see, introduces some additional stability requirements in order for the IMEX schemes to preserve the asymptotic behavior of the equation.
%\begin{remark}~
%\begin{itemize}
%\item Let us emphasize that most of the
%subsequent theory can be generalized to the case where $P(f,f)$ is
%not strictly positive and $\mu$ is an arbitrary nonnegative
%constant. As shown in~\cite{dimarco6} in this case positivity properties can still be achieved for a suitable choice of the time step.
%
%\item The above decomposition correspond to a penalization of the collision term by a BGK operator. In the general case one can introduce a more accurate linear operator to penalize the collision
%operator. Two other possibilities are given by the ES-BGK relaxation
%operator and by the linearized Boltzmann operator.
%\end{itemize}
%\end{remark}

\subsection{Application to the Boltzmann equation}
We can now introduce the general class of IMEX Runge-Kutta schemes for the Boltzmann equation in the form
\begin{eqnarray}%\left\lbrace
%\begin{array}{lll}
&   F^{(i)} = \displaystyle f^{n}+\Delta t \sum_{j=1}^{i-1} \widetilde{a}_{ij} \left(\frac{\mu}{\varepsilon} g(F^{(j)})-v\cdot\nabla_x F^{(j)}\right)+\Delta t\sum_{j=1}^{i} a_{ij}\frac{\mu}{\varepsilon}(M[F^{(j)}]-F^{(j)})\label{eq:GIMEXb} \\
&   f^{n+1} = \displaystyle f^{n}+\Delta t
\sum_{i=1}^{\nu}\widetilde{\omega}_{i}\left(\frac{\mu}{\varepsilon}
g(F^{(i)})-v\cdot\nabla_x F^{(i)}\right)+\Delta
t\sum_{i=1}^{\nu}\omega_{i}\frac{\mu}{\varepsilon}(M[F^{(i)}]-F^{(i)}).
%\end{array}
%\right.
\label{eq:GIMEX1b}
\end{eqnarray}
In the above scheme the explicit method is characterized by the $\nu \times \nu$
matrix $\tA = (\ta_{ij})$, { $\ta_{ij}=0$}, { $j \geq i$} and the
coefficient vectors are { $\tc = (\tc_1,\ldots,\tc_\nu)^T$},
$\tc_i=\sum_{j=1}^{i-1}{\ta_{ij}}$,
$\tw= (\tw_1,\ldots,\tw_\nu)^T$, whereas
the implicit method is a diagonally implicit Runge-Kutta (DIRK)
defined by the $\nu \times
\nu$ matrix $A = (a_{ij})$, $a_{ij}=0$, $j > i$, and the coefficient vectors are { $c =
(c_1,\ldots,c_\nu)^T$}, $c_i=\sum_{j=1}^{\nu}{a_{ij}}$, { $w=
(w_1,\ldots,w_\nu)^T$}.

Let us first recall the definition of asymptotic preserving property~\cite{Jin2}
\begin{definition}
The IMEX scheme (\ref{eq:GIMEXb}-\ref{eq:GIMEX1b}) for the Boltzmann equation is {asymptotic preserving\/} (AP) if in the limit { $\epsilon\to 0$}
the scheme becomes a consistent discretization of the limit system
of the Euler equations (\ref{eq:Euler}).
%We use the notation { $AP_k$} if the scheme
%is of order { $k$} in the limit { $\epsilon\to 0$}.
\end{definition}

Note that if we multiply the IMEX scheme by the vector of collision invariants
$\phi(v)=(1, v, v^2/2)^T$ and integrate in $v$ we get a \emph{moment
scheme} characterized by the explicit method
\begin{eqnarray}
%\left\lbrace
%\begin{array}{lll}
&\displaystyle   \int_{\RR^3} F^{(i)}\phi(v)\,dv = \displaystyle \int_{\RR^3} f^{n}\phi(v)\,dv-\Delta t \sum_{j=1}^{i-1} \widetilde{a}_{ij} \int_{\RR^3} v\cdot \nabla_x F^{(j)}\phi(v)\,dv\label{eq:GIMEXf} \\
&\displaystyle   \int_{\RR^3} f^{n+1}\phi(v)\,dv = \displaystyle
\int_{\RR^3} f^{n}\phi(v)\,dv-\Delta t
\sum_{i=1}^{\nu}\widetilde{\omega}_{i}\int_{\RR^3} v\cdot \nabla_x
F^{(i)}\phi(v)\,dv.
%\end{array}
%\right.
\label{eq:GIMEX1f}
\end{eqnarray}
Thus a sufficient condition for a scheme to satisfy the $AP$
property is that as $\varepsilon\to 0$ we get $F^{(i)}\to
M[F^{(i)}]$, $i=1,\ldots,\nu$ in (\ref{eq:GIMEXb}). In addition we
must require the kinetic numerical solution $f^{n+1}$ to satisfy
some additional numerical stability requirement. We illustrate this
aspect in the sequel.

First let us start with the following Lemma~\cite{DiPa12}
\begin{lemma}
\label{le:apB} If all diagonal element of the triangular coefficient
matrix $A$ that characterize the DIRK scheme in equations
(\ref{eq:GIMEXb}-\ref{eq:GIMEX1b}) are non zero, then \be \lim
_{\varepsilon\rightarrow 0} F^{(i)} =
M[F^{(i)}].\label{eq:IMEXAP2}\ee
\end{lemma}
%\proof
%Consider now the stage $i$ in the original IMEX scheme written in vector form as
%\[
%\varepsilon F=\varepsilon f^n e+\Delta t \tA (\mu g(F)-\varepsilon v\cdot \nabla_x F)+\Delta t A \mu (M[F]-F),
%\]
%with $F=(F^{(1)},\ldots,F^{(\nu)})^T$, $e=(1,\ldots,1)^T$, $g(F)=(g(F^{(1)},\ldots,g(F^{\nu}))^T$ and $M[F]=(M[F^{(1)}],$ $\ldots,M[F^{\nu}])$.
%
%Since $A$ is invertible we can solve the above system for $(M[F^{(i)}]-F^{(i)})$ to get \be \Delta t (M[F^{(i)}]-F^{(i)})=
%\frac{\varepsilon}{\mu}\sum_{j=1}^{i}b_{ij}\left[F^{(j)}-f^{n}+\Delta
%t \sum_{h=1}^{j-1}\widetilde{a}_{jh}\left(v\cdot \nabla_x
%F^{(h)}-\frac{1}{\varepsilon}g(F^{(h)})\right)\right],\ee
%where $A^{-1}=(b_{ij})$. As
%$\varepsilon\rightarrow 0$ we get $F^{(i)}=M[F^{(i)}]$. In fact $\widetilde{A}$ is lower triangular with
%$\widetilde{a}_{ii}=0$ and we have a hierarchy of equations such
%that $F^{(h)}=M[F^{(h)}], h=1,..,j-1$.
%\endproof

Formally Lemma 1 guarantees the AP property of the scheme. However,
as opposite to the case of hyperbolic systems with relaxation, now because of the decomposition of the collision operator the last level (\ref{eq:GIMEX1b}) still depends on
$\varepsilon$. After some manipulations it reads \bea\nonumber
f^{n+1}&=&f^{n}\left(1-\sum_{i,j}{w}_{i}b_{ij}\right)-\Delta
t\sum_{i=1}^{\nu}\widetilde{w}_{i}\left(v\cdot \nabla_x F^{(i)}-\frac{1}{\varepsilon}g(F^{(i)})\right)\\[-.4cm]
\label{eq:fn1}
\\
\nonumber
&+&\Delta t\sum_{i,j,h}{w}_{i}b_{ij}\widetilde{a}_{jh}\left(v\cdot
\nabla_x F^{(h)}-\frac{1}{\varepsilon}g(F^{(h)})\right)
+\sum_{i,j}w_{i}b_{ij}F^{(j)},\eea
where $b_{ij}$ are the elements of $A^{-1}$.
The above expression turns out to be unbounded as $\varepsilon\to 0$ thus originating an unstable scheme.

We introduce the following definition~\cite{BPR11, DiPa12}
\begin{definition}
An IMEX scheme in the form (\ref{eq:GIMEXb})-(\ref{eq:GIMEX1b}) is \emph{globally stiffly accurate} if the following conditions are satisfied
\be
w_{i} = a_{\nu i},\qquad            \widetilde{w}_{i} =
\widetilde{a}_{\nu i},\quad \forall\,i=1,\ldots,\nu.
\label{eq:gsa}
\ee
\end{definition}

We can finally state the main result~\cite{DiPa12}
\begin{theorem}
\label{th:apB} If $ \det A \neq 0$ and the IMEX scheme (\ref{eq:GIMEXb})-(\ref{eq:GIMEX1b}) is globally stiffly accurate, in
the limit $\varepsilon\rightarrow 0$, the IMEX scheme
becomes the explicit RK scheme characterized by
($\widetilde{A}, \widetilde{w}, \widetilde{c}$) applied to the limit
Euler system (\ref{eq:Euler}).
\end{theorem}
%\proof
In order to prove the Theorem it is enough to observe that
the stiffly accurate property implies immediately that
$f^{n+1}=F^{(\nu)}$.

%This can be also seen directly observing that
%as a consequence of (\ref{eq:gsa}) we get
%\[
%\sum_{i,j}{w}_{i}b_{ij}=1,\qquad
%\sum_{i,j,h}{w}_{i}b_{ij}\widetilde{a}_{jh} =
%\widetilde{w}_h,\qquad \sum_{i,j}w_{i}b_{ij}F^{(j)} = F^{(\nu)},
%\]
%which in (\ref{eq:fn1}) yields $f^{n+1}=F^{(\nu)}$ and thus in the limit $\varepsilon\to 0$ we get $f^{n+1}=M[f^{n+1}]$.
%\endproof

\begin{remark}~
\begin{itemize}
\item Theorem above guarantees not only
asymptotic preservation but also asymptotic accuracy, namely the
order of the scheme is preserved in the $\varepsilon\to 0$ limit.
\item
The previous results can be extended to the case of CK-type schemes~\cite{boscarino} with $a_{11}=0$, like the ones considered in~\cite{Filbet}. However in this case asymptotic accuracy holds true only if the initial data are an $O(\varepsilon)$
perturbation of the local Maxwellian equilibrium.
\end{itemize}
\end{remark}

\subsection{Convexity of the schemes}
The determination of general conditions for positivity of
the numerical solution in the space non homogeneous case is quite
difficult. Here we focus on the space homogeneous situation. Not that even in this case due to the reformulation of the collision term the analysis involve the whole IMEX scheme. Moreover the analysis here depends on the particular operator used as a penalization.
%
%In the space homogenous situation the IMEX scheme (\ref{eq:GIMEXb})-(\ref{eq:GIMEX1b}) takes the form
%\begin{eqnarray}
%F^{(i)} &=& \displaystyle f^{n}+\frac{\mu\Delta
%t}{\varepsilon}\sum_{j=1}^{i-1}\widetilde{a}_{ij}\left(\frac{P(F^{(j)},F^{(j)})}{\mu}-M[F^{(j)}]\right)+\frac{\mu\Delta
%t}{\varepsilon}\sum_{j=1}^{i}a_{ij}(M[F^{(j)}]-F^{(j)})
%\\
%f^{(n+1)} &=& \displaystyle f^{n}+\frac{\mu\Delta
%t}{\varepsilon}\sum_{i=1}^{\nu}\widetilde{w}_{i}\frac{P(F^{(i)},F^{(i)})}{\mu}+\frac{\mu
%\Delta t}{\varepsilon}\sum_{i=1}^{\nu}w_{i}(M[F^{(i)}]-F^{(i)})
%\end{eqnarray}

Using the fact that in the space homogenous situation $M[f]$ does not depend on time the IMEX scheme (\ref{eq:GIMEXb})-(\ref{eq:GIMEX1b}) can be rewritten as
\bea
            F^{(i)} &=& \displaystyle
            \sum_{h=1}^{i}\hat{b}_{ih}\left\{\lambda f^n + \sum_{j=1}^{h-1}\widetilde{a}_{hj}\frac{P(F^{(j)},F^{(j)})}{\mu}
            + \displaystyle
            M[f^n]\left(c_h-\widetilde{c}_h\right)\right\}
            \label{eq:IMHO}
            \\
            f^{(n+1)} &=& \displaystyle f^{n}+\frac{\mu\Delta t}{\varepsilon}\sum_{i=1}^{\nu}\widetilde{w}_{i}\frac{P(F^{(i)},F^{(i)})}{\mu}
            +\frac{\mu\Delta t}{\varepsilon}\sum_{i=1}^{\nu}w_{i}(M[f^n]-F^{(i)})
            \label{eq:IMHOb}
          \eea
where $\lambda=\varepsilon/(\mu\Delta t)$ and $\hat{b}_{ij}$ are the elements of $(\lambda I + A)^{-1}$.

The following theorem gives sufficient conditions for the above expression to represent a convex combination of probability densities.

\begin{theorem}
A sufficient condition to guarantee that $f^{n+1}\geq 0$ when $f^n\geq 0$ in (\ref{eq:IMHO})-(\ref{eq:IMHOb}) is that the scheme is
globally stiffly accurate and the following conditions holds true
\be 0\leq\sum_{h=1}^{i}\hat{b}_{ih}c_h\leq 1,\quad 0\leq
\sum_{h=1}^{i}\hat{b}_{ih}\left(c_{h}-\widetilde{c}_h\right)\leq 1,\quad\forall\,\,
i=1,\ldots,\nu.\ee
\be 0\leq \sum_{h=j+1}^{i}\hat{b}_{ih}\widetilde{a}_{hj}\leq 1,\quad \forall\,\,
i=1,\ldots,\nu,\quad j=1,\ldots,i-1.\ee
\end{theorem}
Since the result is based on a convexity argument we also have an entropic result for the schemes.
\begin{corollary}
Under the assumptions of Theorem 1, if in addition the operator $P(f,f)$ satisfies
\be
H\left(\frac{P(f,f)}{\mu}\right)\leq H(f),\qquad H(f)=\int_{\R^3} f\log f\,dv,
\ee
then $H(f^{n+1})\leq H(f^n)$.
\end{corollary}
%
%
%\begin{remark}~
%\begin{itemize}
%\item Since the result is based on a contractivity argument it follows that the IMEX schemes are also entropic provided we have an
%estimate of the type~\cite{Vi}
%\be
%H(P(f,f))\leq H(f),\qquad H(f)=\int_{\R^3} f\log f\,dv.
%\ee
%
%
%\item As a consequence this theorem gives sufficient condition for a scheme based on a splitting method in order to be fully positive. This permits, for example, the use of Monte Carlo methods in combination with IMEX methods~\cite{dimarco6}.
%
%\end{itemize}
%\end{remark}
\begin{figure}
\begin{center}
\includegraphics[scale=0.4]{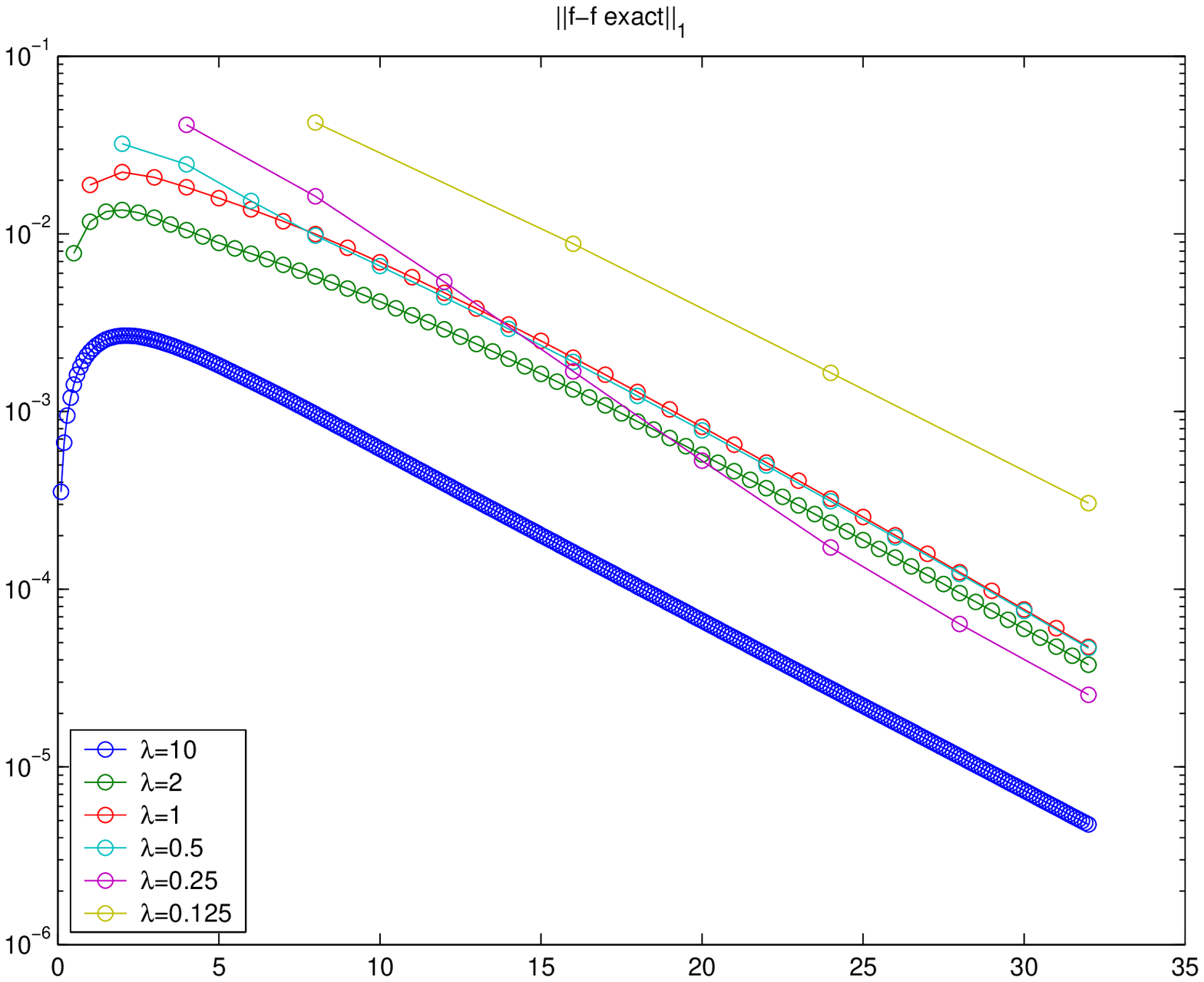}
\includegraphics[scale=0.4]{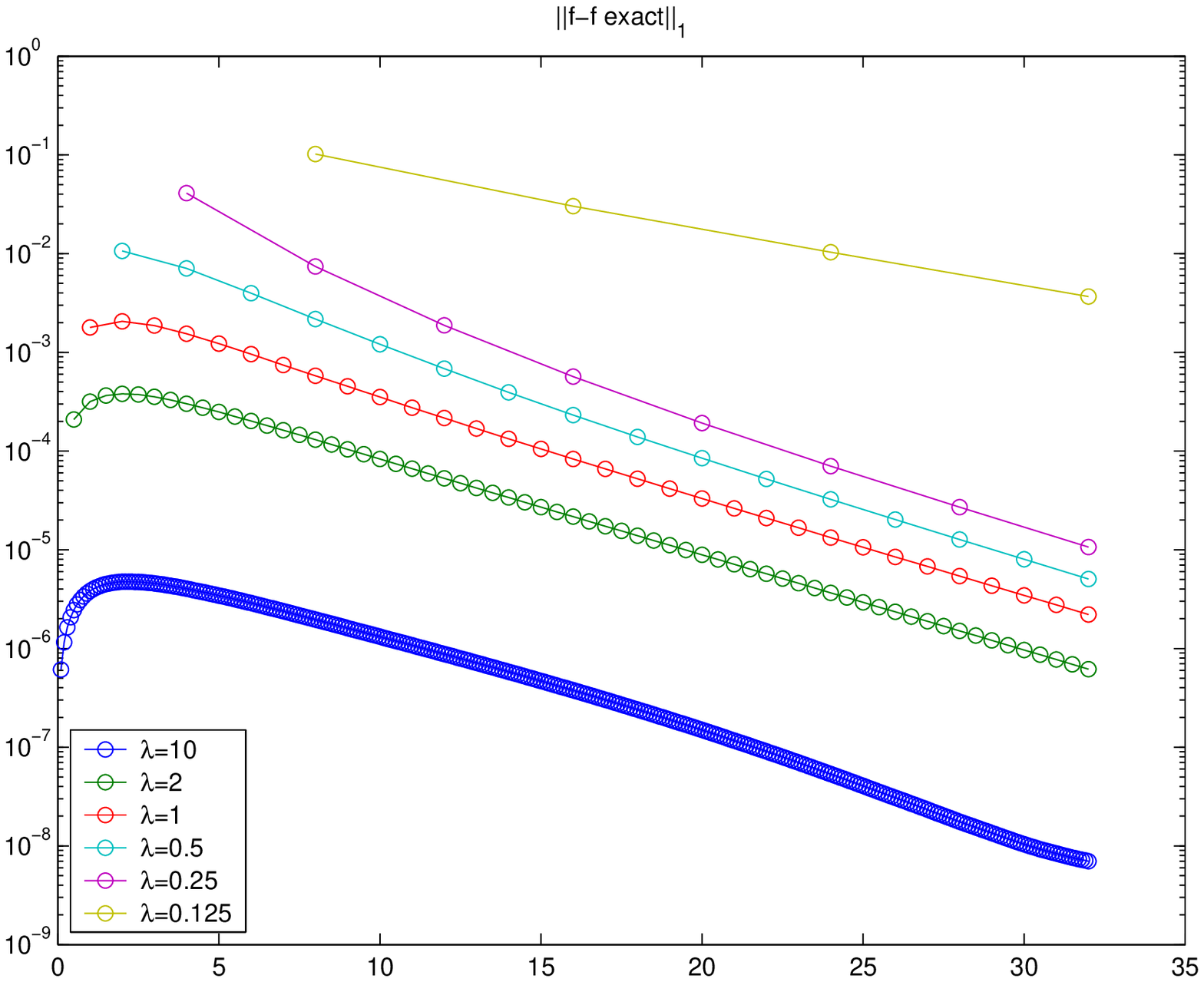}
\caption{$L_{1}$ error for the distribution function $f$ for the
second (left) and the third (right) IMEX-BE
method.}\label{fig:error1}
\end{center}
\end{figure}
\section{Examples and numerical results}
In Table 1 we report one example of a second order
asymptotic preserving scheme~\cite{DiPa12}. The scheme is also positivity
preserving for $\lambda \leq 1$ and asymptotically
accurate. We also report in Table 2 a third
order scheme globally stiffly accurate~\cite{boscarino}.

The schemes can be schematically summarized using a \emph{double
Butcher tableau} of the type~\cite{PRimex}
\[
    \begin{array}{c|cccc}
    \tc & \tA \\
    \hline\\[-.2cm]
    & \tw^T
    \end{array}\qquad\qquad
\begin{array}{c|cccc}
    c & A \\
    \hline\\[-.2cm]
    & w^T
    \end{array}
\]
Note that although the schemes use several implicit evaluations they
are still optimal in terms of number of evaluation of the collision
operator. This, in fact, is characterized only by the number of
explicit function evaluations. We used the notation
name$(k,\sigma_E,\sigma_I)$ where $k$ is the order and $\sigma_E,\sigma_I$ characterize the number of
evaluations of the explicit and implicit schemes respectively.

The numerical test is an homogeneous relaxation problem in the two
dimensional velocity space. The molecules are Maxwellian and a
fast spectral method~\cite{MP} is used to compute the collision operator with
$N_{v}=64$ grid points in each velocity direction and a grid
$[-v_{\max}, v_{\max}]^2$ with $v_{\max}=3\pi$.
%The non equilibrium
%initial data is given by
%\[
%f(v,0)=\frac{v^{2}}{\pi
%\sigma^{4}}\exp\left(-\frac{v^{2}}{\sigma^{2}}\right),
%\]
In this case, the exact solution is given
by
\[
f(v,t)=\frac{1}{2\pi
S^{2}\sigma^{2}}\left(2S-1+\frac{1-S}{2S}\frac{v^{2}}{\sigma^2}\right)\exp\left(-\frac{v^{2}}{2S\sigma^{2}}\right),
\qquad S(t)=1-\frac{\exp\left(-\sigma^{2}t/8\right)}{2},
\]
where we took $\sigma=1$.
The figure \ref{fig:error1} shows the error of the schemes for different choices of the time step $\Delta t$ (stability condition for the explicit Euler scheme is $\Delta t=1$). We can clearly observe the expected accuracy of the schemes even for large time steps.
% and the asymptotic preserving property since the error vanishes for %larger times.
%On the right we reported the error for the fourth order moment
%(exact and computed solution) again for different values of the time
%step. For this test case, the stability condition for the explicit
%Euler scheme imposes $\Delta t=1$.
%We clearly see that the IMEX
%scheme is stable and asymptotic preserving, the error, in fact, goes
%to zero in the limit $t\rightarrow \infty$.

%\begin{table}[t]
%{\small
%\[
%    \begin{array}{c|cccc}
%      0 & 0        & 0              \\
%      1   & 1 & 0 \\
%      \hline
%          & 1   &  0
%    \end{array}\qquad
%\begin{array}{c|cccc}
%      1 & 1        & 0              \\
%      1   & 0 & 1 \\
%      \hline
%          & 0   &  1
%    \end{array}
%\]
%\caption{Tableau of the first order IMEX-BE1(2,2,1) asymptotic and positivity
%preserving IMEX scheme}
%\label{tb:firstB} }
%\end{table}

\begin{table}[t]
{\small
\[
    \begin{array}{c|ccccc}
      0 & 0 & 0 & 0 & 0   \\
      0   & 0 & 0 & 0 & 0 \\
      1   & 0 & 1 & 0 & 0 \\
      1   & 0 & 1/2 & 1/2 & 0 \\
      \hline
          & 0   &  1/2 & 1/2 & 0
    \end{array}\qquad
  \begin{array}{c|ccccc}
      2 & 2 & 0 & 0 & 0   \\
      0   & -2 & 2 & 0 & 0 \\
      1   & 0 & -1 & 2 & 0 \\
      1   & 0 & 1/2 & -3/2 & 2 \\
      \hline
          & 0   &  1/2 & -3/2 & 2
    \end{array}\]
\caption{Tableau of the second order IMEX-BE(2,2,4) asymptotic and positivity
preserving IMEX scheme.}
\label{tb:secondB} }
\end{table}

\begin{table}[t]
{\small
\[
    \begin{array}{c|cccccc}
      0 & 0 & 0 & 0 & 0 & 0  \\
      1   & 1 & 0 & 0 & 0 & 0\\
      2/3   & 4/9 & 2/9 & 0 & 0 & 0\\
      1   & 1/4  & 0 &  3/4 & 0 & 0\\
      1   & 1/4  & 0 &  3/4 & 0 & 0\\
      \hline
          & 1/4  & 0 &  3/4 & 0 & 0
    \end{array}\qquad
  \begin{array}{c|cccccc}
      0 & 0 & 0 & 0 & 0 & 0 \\
      1 & 1/2 & 1/2 & 0& 0 & 0   \\
      2/3 & 5/18 & -1/9 & 1/2 & 0 & 0 \\
      1   & 1/2 & 0 & 0 & 1/2 & 0 \\
      1 & 1/4 & 0  &  3/4 & -1/2 & 1/2\\
      \hline
          & 1/4 & 0  &  3/4 & -1/2 & 1/2
    \end{array}\]
\caption{Tableau of the third order IMEX-BE(3,5,5) globally stiffly accurate IMEX scheme.}
\label{tb:thirddB} }
\end{table}
\vskip 1cm \textbf{Acknowledgement}. G. Dimarco was supported by the
French ANR project BOOST.
%The authors would like to thank Professor
%S. Jin and Prof. F. Filbet for the stimulating discussions.
%%%%%%%%%%%%%%%%%%%%%%%%%%%%%%%%%%%%%%%%%%%%%%%%%%%%%%%%%%%%%%%%%%%%%%%%%%%%%%%%%%%%%%%%%%%%%%%%%%%%%%%%%%%%%%%%%

\end{document}